\input ./style/arxiv-general.cfg
\documentclass[aos,MSNbibl,dvips]{arximspdf}
\makeatletter
   \@ifpackageloaded{graphicx}{}{\usepackage{graphicx}}
\makeatother

\doi{10.1214/14-AOS1270A}
\volume{43}
\issue{4}
\pubyear{2015}
\firstpage{1429}
\lastpage{1436}
\docsubty{FLA}
\referstodoi{10.1214/14-AOS1270}

\makeatletter
\newcommand{\rrvert}{\vert}
\newcommand{\rrVert}{\Vert}
\newcommand{\llvert}{\vert}
\newcommand{\llVert}{\Vert}
\makeatother

\begin{document}
\begin{frontmatter}
\vspace*{12pt}
\title{Discussion of ``Frequentist coverage of adaptive nonparametric
Bayesian credible sets''}
\runtitle{Discussion}

\begin{aug}
\author[A]{\fnms{Richard}~\snm{Nickl}\corref{}\ead[label=e1]{r.nickl@statslab.cam.ac.uk}}
\runauthor{R. Nickl}
\affiliation{University of Cambridge}
\address[A]{Statistical Laboratory\\
Department of Pure Mathematics\\
\quad and Mathematical Statistics\\
University of Cambridge\\
CB3 0WB Cambridge\\
United Kingdom\\
\printead{e1}}
\end{aug}

%
\received{\smonth{9} \syear{2014}}


\end{frontmatter}

\section{Introduction}

I would like to congratulate Botond Szab\'o, Aad van der Vaart and
Harry van Zanten \cite{SVV13} for a fundamental and thought provoking
article on a highly important topic. One of the key contributions of
statistics to modern science may arguably be \textit{the theory of
uncertainty quantification}. Assessing the accuracy of an estimate by a
confidence statement goes beyond the mere search for an efficient
statistical algorithm. In particular, within the contemporary search
for \textit{adaptive} procedures, research of the last decade has
revealed that the construction of adaptive confidence statements is
fundamentally \textit{harder}---in an information theoretic sense---than the construction of adaptive algorithms. Confidence statements are
at the same time crucial for the main application of modern data
analysis, which is to accept or reject hypotheses.

Szab\'o, van der Vaart and van Zanten tackle the important topic as to
whether increasingly popular Bayesian methodology can actually provide
objective uncertainty quantification methods in nonparametric models or
not. The nonparametric setting is a key test-case for the general
paradigm of high-dimensional modeling that has emerged recently in statistics.

My discussion of the paper surrounds the two focal points of why
``Bayesian uncertainty quantification'' is a mathematically and
conceptually nontrivial problem: the first has nothing to do with
adaptation and addresses some of the probabilistic subtleties intrinsic
to the Bayesian approach to provide ``credible sets.'' The second point
is common to all frequentist procedures and is about the fact that
adaptive uncertainty quantification is in general only possible under
``signal-strength'' conditions on the underlying parameter.

\section{Freedman's paradox and the nonparametric Bernstein--von Mises theorem}

I first want to discuss the fact that the frequentist coverage
probabilities obtained by Szab\'o, van der Vaart and van Zanten for
their credible sets are not \textit{exact}, that is to say, not of the
precise asymptotic level $1-\alpha$, and the related question of why
obtaining exact posterior asymptotics in the nonparametric situation is
a subtle matter.

Consider observations $Y \sim P_\theta$ with parameter space $\theta
\in\Theta$, a prior $\Pi$ on $\Theta$ and resulting posterior
distribution $\Pi(\cdot|Y)$ of $\theta|Y$. The classical
finite-dimensional ($\Theta\subseteq\mathbb R^p$) \textit{Bernstein--von Mises theorem} asserts---under fairly mild assumptions on $\Pi$
and on the parameterization $\theta\mapsto P_\theta$---that we have
approximately
\[
\mathcal L\bigl(\sqrt n (\theta- \bar\theta)|Y\bigr) \approx N\bigl(0, I(
\theta_0)^{-1}\bigr) \qquad\mbox{when } Y \sim
P_{\theta_0}, \theta_0 \in\Theta.
\]
Here $\bar\theta= \bar\theta(Y)$ is any efficient estimator of
$\theta$ such as the maximum likelihood estimator (MLE) or the
posterior mean $E(\theta|Y)$, $I(\theta_0)$ is the Fisher
information, and the approximation holds in the small noise or large
sample limit, in total variation distance. As a consequence, computing
posterior probabilities is approximately equivalent to computing
``optimal frequentist'' probabilities under $N(\bar\theta, I(\theta
_0)^{-1}/n)$, and the natural level $1-\alpha$ Bayesian credible set
%
\begin{equation}
\label{parcred} C_n = \bigl\{\theta\dvtx  \bigl\llVert \theta- E(\theta|Y)
\bigr\rrVert \le r_{\alpha,n}\bigr\}\qquad \mbox{with } r_{\alpha,n}
\mbox{ s.t. }\Pi(C_n|Y) = 1-\alpha
\end{equation}
asymptotically coincides with the classical one based on the MLE. In
particular, we have frequentist coverage
%
\begin{equation}
\label{parcov} P_{\theta_0}(\theta_0 \in C_n) \to1-
\alpha\qquad\mbox{as } n \to \infty.
\end{equation}
For the frequentist the main idea behind this phenomenon is similar to
the bootstrap: if $Y \sim P_{\theta_0}$, the (known) posterior
distribution of $\theta|Y - E(\theta|Y)$ serves as a proxy for the
(unknown) distribution of $E(\theta|Y) - \theta_0$.

In his influential 1999 Wald lecture, Freedman \cite{Free99} has shown
that in the case where $\Theta$ is infinite-dimensional, the above
phenomenon need not occur. Freedman considered precisely the setting of
Szab\'o, van der Vaart and van Zanten: in the standard nonparametric
sequence space model
%
\begin{equation}
Y_k = \theta_k + \frac{1}{\sqrt n} g_k,\qquad k
\in\mathbb N;\qquad g_k \stackrel{\mathrm{i.i.d.}}{\sim} N(0,1),\qquad  \theta\in
\ell_2,
\end{equation}
one considers natural conjugate Gaussian priors
%
\begin{equation}
\label{gp} \Pi= \bigotimes_{k \in\mathbb N} N\bigl(0,
k^{-1-2\gamma}\bigr),\qquad \gamma>0,
\end{equation}
for a $\gamma$-regular signal $\theta$. Freedman then showed that
even when the true signal $\theta_0$ is $\beta$-regular with $\beta
>\gamma$---so in a favorably ``well-specified'' situation where
$\beta$ is known---the natural posterior credible set paralleling
(\ref{parcred}),
%
\begin{eqnarray}\label{parcred2}
C_n = \bigl\{\theta\in\ell_2\dvtx  \bigl
\llVert \theta- E(\theta|Y)\bigr\rrVert _{\ell_2} \le r_{\alpha, n}\bigr
\}
\nonumber\\[-8pt]\\[-8pt]
\eqntext{\mbox{with } r_{\alpha, n}\mbox{ s.t. }\Pi(C_n|Y) =
1-\alpha,}
\end{eqnarray}
does in fact \textit{not} satisfy (\ref{parcov}) as $n \to\infty$,
rather the frequentist and Bayesian variances of $\llVert  \theta- E(\theta
|Y)\rrVert  ^2_{\ell_2}$ scale differently and the Bernstein--von Mises
theorem does not hold in this infinite-dimensional setting. See
Freedman's original paper \cite{Free99} and also Leahu \cite{leahu11}
for a recent account.

We note that this ``paradox'' has nothing to do with ``adaptation''
(since $\beta$ is known above), but is a mathematical artefact of the
Bayesian formalism to construct credible sets in the
infinite-dimensional setting. It prevents Bayesian $1-\alpha$ credible
balls in $\ell_2$ from being asymptotically \textit{exact}
frequentist $1-\alpha$ confidence balls. The pathology occurs because
one insists on ``exact'' asymptotic level $1-\alpha$, and the results
by Szab\'o, van der Vaart and van Zanten show that if one ``blows up
the Bayesian radius $r_{\alpha, n}$'' by a factor of $L$ as in
equation (3.2) of their paper \cite{SVV13}, then as $L \to\infty$
this allows to obtain ``conservative'' frequentist confidence sets in
the sense that, as $n \to\infty$,
\[
P_{\theta_0}(\theta_0 \in C_n) \to1 \ge1-\alpha.
\]
While such a construction is theoretically satisfactory from a
frequentist point of view, this approach has a practical drawback: in
applications one does not know how to choose $L$ and prefers to have a
simple, fully Bayesian, rule that discards 5\% of all posterior draws
and uses the remaining 95\% to graphically describe a credible region.

In the recent papers \cite{CN13,CN14} by Castillo and myself, a new
approach to nonparametric Bernstein--von Mises theorems has been put
forward. In essence, the idea is to modify the geometry of the credible
set in (\ref{parcred2}) and to replace $\ell_2$-balls by other
shapes. These shapes correspond to norms in sequence space that induce
weaker topologies than $\ell_2$ and for which a ``weak functional
Bernstein--von Mises theorem'' can be proved. For instance, if we
consider ellipsoids in a sequence space of the form
%
\begin{equation}
\label{ellips} \mathcal E(M) = \biggl\{(\theta_k)\dvtx  \sum
_{k} \frac{\theta_k^2}{w_k} \le M^2 \biggr\},\qquad
\frac{w_k}{k (\log k)^{\delta}} \uparrow\infty, \delta>1, 0<M<\infty
\end{equation}
or, in case the sequence space model corresponds to a double-indexed
basis $\{e_{lk}\dvtx  l \in\mathbb N \cup\{0\}, k =0, \dots, 2^l-1\}$,
multi-scale sup-norm balls of the form
%
\begin{equation}
\mathcal E(M) = \biggl\{(\theta_{lk})\dvtx  \sup_{l \ge0}
\frac{\max_k
\llvert \theta_{lk}\rrvert }{w_l} \le M \biggr\},\qquad \frac{w_l}{\sqrt l} \uparrow \infty, 0<M<\infty,
\end{equation}
then, under mild assumptions on the prior, \cite{CN13,CN14} prove
that, as $n \to\infty$,
%
\begin{equation}
\mathcal L\bigl(\sqrt n (\theta- \bar\theta)|Y\bigr) \to\mathcal N\qquad
\mbox{weakly in } \mathbb H,
\end{equation}
in $P_{\theta_0}$-probability. Here $\mathbb H$ are the sequence
spaces that have norm balls $\{\theta\dvtx  \llVert  \theta\rrVert  _\mathbb H \le1\} =
\mathcal E(1)$, $\mathcal N$ is the Gaussian measure on $\mathbb H$
corresponding to a pure white noise $\bigotimes_{k \in\mathbb N}
N(0,1)$, and $\bar\theta= \bar\theta(Y)$ is equal to the maximum
likelihood estimator $Y=(Y_k\dvtx  k \in\mathbb N)$ or to the posterior
mean $E(\theta|Y)$.

As a consequence of weak convergence toward $\mathcal N$, one can show that
%
\begin{equation}
\sup_M \bigl\llvert \Pr\bigl(\sqrt n \bigl(\theta-E(
\theta|Y)\bigr) \in\mathcal E(M) |Y\bigr) - \mathcal N\bigl(\mathcal E(M)\bigr)
\bigr\rrvert \stackrel{P_{\theta_0}}{\to} 0
\end{equation}
as $n \to\infty$, and from this \cite{CN13,CN14} deduce that
credible sets
\[
C_n = \bigl\{\theta\dvtx  \bigl\llVert \theta-E(\theta|Y)\bigr\rrVert
_{\mathbb H} \le r_{\alpha,
n}\bigr\}\qquad\mbox{where } r_{\alpha, n}
\mbox{ is such that } \Pi(C_n|Y)=1-\alpha %
\]
have correct asymptotic frequentist coverage: as $n\to\infty$,
\[
P_{\theta_0}(\theta_0 \in C_n) \to1-\alpha.
\]

Two main questions arise from this construction, one theoretical, one
practical. The theoretical one asks whether such confidence sets can
reconstruct nonparametric signals in a minimax optimal way. In \cite
{CN13,CN14} it is shown that this can be the case by using
high-frequency information in the posterior appropriately. This can be
extended to the adaptive setting (see Ray \cite{R14}), where
nonparametric Bernstein--von Mises theorems in $\mathbb H$ are proved
for the empirical Bayes procedure of Szab\'o, van der Vaart and van Zanten.

The second question is as follows: is such a construction practical,
and do such credible sets look substantially different from the
(perhaps) more intuitive $\ell_2$-type credible sets? From a
computational point of view the sets $C_n$ are quite tractable: for
instance, in the multi-scale case the computation of $C_n$ consists of
finding constants $r_{\alpha, n}$ such that
\[
\frac{|\theta_{lk}- E(\theta_{lk}|Y)|}{w_l} \le r_{\alpha,n} \qquad \forall k,l
\]
happens for $(1-\alpha) \times100\%$ of the posterior draws. The
theory in \cite{CN14} implies
\[
\sqrt n \cdot r_{\alpha,n} \stackrel{P_{\theta_0}}{\to} \operatorname{const} \neq0,
\]
and so a multi-scale posterior credible ball has a natural
interpretation as a simultaneous credible set for a large class of
semi-parametric coordinate projection functionals obtained from
thresholding each projection at a level slightly larger than $1/\sqrt
n$ (recalling that $w_l$ is slightly larger than $\sqrt l$).

More concretely, simulations by Ray \cite{R14} show that the
differences to the standard $\ell_2$-approach are marginal in several
practical examples; see Figure~\ref{fig1} below.

\begin{figure}

\includegraphics{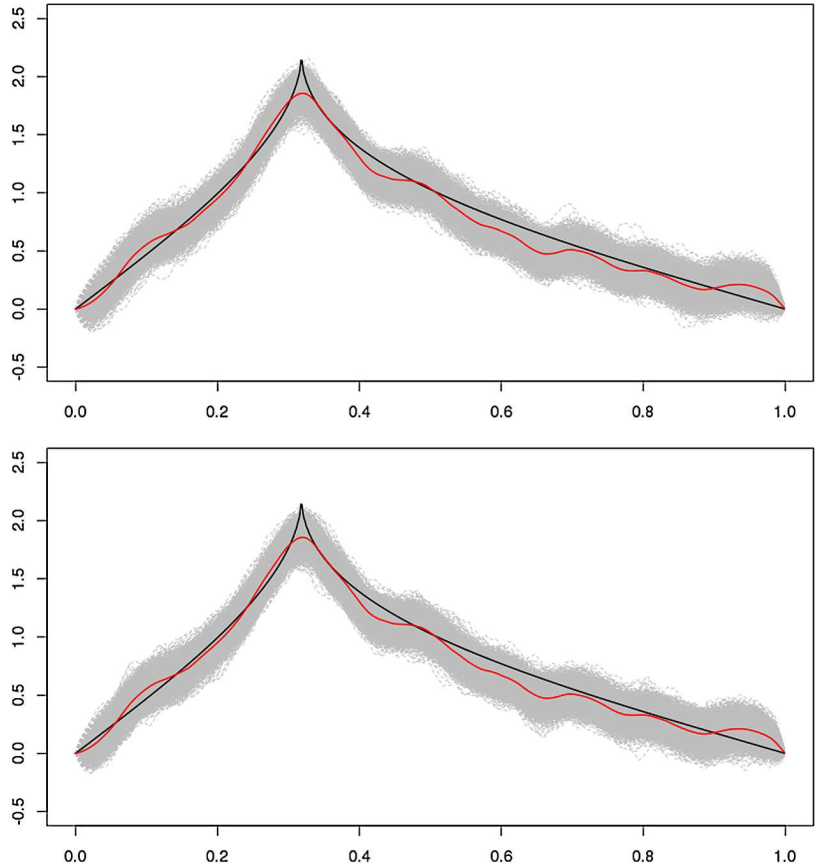}

\caption{The top display shows a credible set generated from
ellipsoids (\protect\ref{ellips})
and the bottom from an $\ell_2$ ball. Both credible sets are
based on observations in Gaussian white noise ($n=1000$), and with
Gaussian priors, with 100,000 posterior draws plotted as grey clouds.
The red curve depicts the posterior mean and the black curve the true
function. See \cite{R14} for details and more extensive simulation results.}\label{fig1}
\end{figure}

It is striking to observe that, although the norms $\llVert  \cdot\rrVert  _{\ell
_2}$ and $\llVert  \cdot\rrVert  _\mathbb H$, as well as the rules these norms
induce to accept or reject posterior draws in the construction of a
credible set, are quite different, the visualized credible sets of both
approaches look very similar.

It is also worthwhile noting that in both cases the credible ball
actually ``covers the true function'' despite the graph suggesting that
pointwise coverage fails. The reason is that $\ell_2$-confidence balls
are insensitive to lack of coverage on intervals of small Lebesgue
measure. A frequentist theory for simultaneous credible ``bands'' is
thus also of interest---some first results in this direction are given
in \cite{CN14,R14}.

I am unsure to which extent $\ell_2$-credible sets are ``applied in
current practice'' as claimed in the Introduction of \cite{SVV13},
particularly if one has to choose ``blow-up'' constants $L$.
Practitioners may prefer to avoid such choices, and instead compute
posterior credible balls in $\mathbb H$-spaces. At any rate, it remains
a mathematical fact that the nonparametric Bernstein--von Mises theorem
\textit{does} hold in the spaces $\mathbb H$, whereas it does \textit
{not} hold in $\ell_2$.

\section{``Honest'' nonparametric models and ``self-similar'' functions}

Perhaps more important than the question of how to obtain \textit
{exact} coverage statements for Bayesian credible sets (discussed in
the previous section) is the question of existence of \textit
{adaptive} confidence sets---Bayesian or not. It is one of the more
surprising insights of the theory of nonparametric and high-dimensional
inference that estimators that adapt to unknown regularity properties
(such as smoothness or sparsity) exist, whereas associated confidence
sets in general \textit{do not}. Roughly speaking, the reason behind
this is that an ``honest'' ($=$uniform in the parameter $\theta$)
adaptive confidence set implicitly solves the testing problem of
whether a signal belongs to a given regularity class or not, and that
such tests simply do not exist over the entire parameter spaces
considered in nonparametric \textit{estimation}. Rather, some kind of
\textit{signal-strength condition} needs to be enforced on the
elements of the parameter space to construct confidence sets for
adaptive estimators. See Hoffmann and Nickl \cite{HN11} and Nickl and
van de Geer \cite{NvG13} for two basic instances of this fact (in
nonparametrics and sparse regression, resp.).

Among such signal-strength conditions, the ``self-similarity''
assumptions introduced in Gin\'e and Nickl \cite{GN10} have proved
compatible with commonly used adaptive frequentist procedures (such as
Lepski's method). In the $L^\infty$-setting of confidence bands they
are also shown to be necessary (see \cite{B12}) if one wants to adapt
to a continuum of smoothness parameters, as is usually the case in
nonparametric statistics. The starting point of Szab\'o, van der Vaart
and van Zanten is to transpose the $L^\infty$-type self-similarity
condition from \cite{GN10} into their $\ell_2$-risk setting:
%
\begin{equation}
\label{ssim} \sum_{k=N}^{\rho N}
\theta_k^2 \ge\varepsilon\llVert \theta\rrVert
^2_{S^\beta
}N^{-2\beta}\qquad \forall N \ge N_0
\mbox{ with ``tolerance'' factor } \varepsilon>0,
\end{equation}
whenever $\theta$ belongs to a Sobolev space $S^\beta$ with norm
\[
\llVert \theta\rrVert ^2_{S^\beta} = \sum
_k \theta_k^2 k^{2\beta},\qquad S^\beta
(B) = \bigl\{\theta\dvtx  \llVert \theta\rrVert _{S^\beta} \le B\bigr\}.
\]
Here $\rho>2, N_0 \in\mathbb N$ are fixed constants; see
equation~(3.4) in \cite{SVV13}. Note that finiteness of the Sobolev
norm implies
%
\begin{equation}
\label{upbd} \sum_{k \ge N} \theta_k^2
\le\llVert \theta\rrVert _{S^\beta}^2 N^{-2\beta
}\qquad
\forall N \in\mathbb N,
\end{equation}
and the idea behind (\ref{ssim}) is hence that over repeated blocks $\{
N, \dots, \rho N\}$ the signal $\theta$ indicates that it is actually
exactly $\beta$-regular. A nice observation of Szab\'o, van der Vaart
and van Zanten is that this condition can in fact be substituted by the
slightly more general ``polished tail'' condition
%
\begin{equation}
\label{pt} \sum_{k=N}^{\rho N}
\theta_k^2 \ge L_0^{-1}\sum
_{k \ge N} \theta _k^2 \qquad \forall
N \ge N_0\mbox{ for some } L_0>0,
\end{equation}
which effectively means that the blocks in (\ref{ssim}) have, for
every $N$ large enough and up to a small constant $L_0^{-1}$, the same
signal strength as the full tail series $\sum_{k \ge N} \theta_k^2$.
This condition is conceptually somewhat cleaner than (\ref{ssim}), as
it does not require the identification of the unknown regularity
parameter $\beta$, although it implicitly does so in the sense that
(\ref{pt}) implies that (\ref{ssim}) and (\ref{upbd}) cannot hold
for multiple values of $\beta$.

\begin{quote}
The key issue I want to discuss here is in which sense exactly
conditions like (\ref{ssim}) or~(\ref{pt}) are necessary for adaptive
inference procedures to exist in the setting of the paper \cite{SVV13}
under discussion.
\end{quote}

Since Szab\'o, van der Vaart and van Zanten are considering $\ell
_2$-risk, the situation is qualitatively different from the $L^\infty
$-setting for which the lower bounds in \cite{B12} apply. First of
all, as also noted by the authors, when adaptation is sought after for
$\beta$ contained in fixed smoothness windows $[\beta_0, 2\beta_0]$,
a direct construction of an adaptive confidence set is possible without
any restrictions on the parameter space. However, the constraint $\beta
\in[\beta_0, 2 \beta_0]$ is not satisfactory in the typical
situations of nonparametric inference. Once relaxed,
information-theoretic arguments imply that restrictions on the
parameter space $S^{\beta}$ become necessary (e.g., Theorems 1 or 4 in
\cite{BN13}). Employing conditions of the kind (\ref{ssim}) or (\ref
{pt}) to enforce such restrictions, one notices that these assumptions
can be weakened quantitatively by increasing the windows $[N, \rho N]$
over which the lower bound of the signal is allowed to accrue. The
question arises whether the window sizes $[N, \rho N]$ with $\rho>2$
are \textit{minimal} conditions for the existence of adaptive
confidence sets or whether larger windows are admissible, pertaining to
larger parameter spaces for which inference is possible. \textit{For
self-similar classes it is shown in} \cite{NS14} \textit{that condition}
(\ref{ssim}) \textit{is} \textit{not} \textit{optimal}, \textit{and that in turn}
(\ref{pt}) \textit{can also not be}.

Let us describe the results from \cite{NS14} to understand in what
sense weaker conditions are possible: let $N_0 \in\mathbb N,
0<b<B<\infty$. For $\varepsilon\in(0,1]$ and $c_\beta=16\times
2^{2\beta+1}$, define the set
%
\begin{equation}
\label{ssim2} S^{\beta}_{\varepsilon} = \Biggl\{\theta\in
S^\beta(B)\dvtx  \sum_{k {=}
N^{(1-\varepsilon)}}^{N}
\theta_{k}^2 \ge c_\beta\llVert \theta\rrVert
_{S^\beta
}^2 N^{-2\beta}\ \forall N \ge N_0
\Biggr\}.
\end{equation}
Again, as in (\ref{ssim}), sufficiently large signal blocks have to
appear repeatedly. But now these blocks are allowed to have increased
\textit{window-width} since
\[
N^{\varepsilon} \gg\rho\qquad\mbox{as } N \to\infty,
\]
and allow for an \textit{asymptotically shrinking tolerance factor}
\[
\varepsilon=N^{-2\varepsilon\beta} c_\beta\to0 \qquad\mbox{as } N \to\infty
\]
in the lower bound. \textit{In particular}, (\ref{ssim2}) \textit{only
approximately identifies the smoothness of} $\theta$ \textit{in the sense that
it can be satisfied}, \textit{unlike} (\ref{ssim}) \textit{or} (\ref{pt}), \textit{for multiple
values of} $\beta$ \textit{simultaneously}.

As shown in \cite{NS14}, signal strength conditions enforced through
(\ref{ssim2}) allow for the construction of honest adaptive confidence
$\ell_2$-balls for signals
\[
\theta\in\bigcup_{\beta_{\min} \le\beta\le\beta_{\max}} S^\beta_{\varepsilon(\beta)},\qquad 0<
\beta_{\min}<\beta_{\max}<\infty,
\]
under (effectively) the following conditions on $\varepsilon$:
\[
\varepsilon(\beta) <\tfrac{1}{2}\qquad \forall\beta\in[\beta_{\min
},
\beta_{\max}]\mbox{ is necessary,}
\]
whereas
\[
\varepsilon(\beta) < \frac{\beta}{2\beta+1/2} \qquad  \forall\beta\in [\beta_{\min},
\beta_{\max}] \mbox{ is sufficient}.
\]
Note that $\beta_{\min}<\beta_{\max}$ are arbitrary, and hence the
lower bound cannot be improved in general, since in the limit $\beta
\to\infty$ we have $\beta/(2\beta+1/2) \to1/2$.

We conclude that requiring lower bounds in windows of size $[N, \rho
N]$ as in (\ref{ssim}), (\ref{pt}) is too strong a requirement for
adaptive $\ell_2$-confidence sets, and the results in the paper by
Szab\'o, van der Vaart and van Zanten are suboptimal from an
information-theoretic perspective. It would be interesting to know
whether this suboptimality is an artefact of the proofs or of the
particular Bayesian inference procedure used, although it may be
difficult to answer this question.


\section*{Acknowledgment}
I would like to thank Kolyan Ray for allowing
me to reproduce Figure~1 from his paper \cite{R14}.


%

\printaddresses
\end{document}